\documentclass[reqno,11pt]{amsart}
\usepackage{amsmath}
\usepackage{latexsym}
\usepackage{amsfonts}
\newtheorem{thm}{Theorem}[section]
\newtheorem{lem}[thm]{{Lemma}}

\newtheorem{prop}[thm]{{Proposition}}

\theoremstyle{definition}

\newtheorem{ex}[thm]{{\bf Example}}
 \numberwithin{equation}{section}
\begin{document}
\newcommand{\num}{S}
\newcommand{\den}{B}
\newcommand{\hTheta}{\widehat{\Theta}}
\newcommand{\cS}{{\mathcal {S}}}
\newcommand{\cB}{{\mathcal {B}}}
\newcommand{\cE}{{\mathcal {E}}}
\newcommand{\PR}{{\bf IP}_{\kappa}}
\renewcommand{\theequation}{\thesection.\arabic{equation}}
\newcommand{\C}{{\mathbb C}}
\newcommand{\D}{{\mathbb D}}
\newcommand{\T}{{\mathbb T}}
\newcommand{\N}{{\mathbb N}}

\title[Multi-point interpolation]{On a multi-point 
interpolation problem for generalized Schur functions}
\author{Vladimir Bolotnikov}
\address{Department of Mathematics\\
The College of William and Mary \\
Williamsburg, VA 23187-8795, USA}
\email{vladi@math.wm.edu}   

\begin{abstract}
The nondegenerate Nevanlinna-Pick-Carath\'eodory-Fejer interpolation 
problem with finitely many interpolation conditions always has 
infinitely many solutions in a generalized Schur class $\cS_\kappa$ for 
every $\kappa\ge \kappa_{\rm min}$ where the integer $\kappa_{\rm min}$ 
equals the number of negative eigenvalues of the Pick matrix associated 
to the problem and completely determined by interpolation data. A linear 
fractional description of all $\cS_{\kappa_{\rm min}}$ solutions of the 
(nondegenerate) problem is well known. In this paper, we present a 
similar result for an arbitrary $\kappa\ge \kappa_{\rm min}$.
\end{abstract}
\maketitle

\centerline{\em Dedicated to Professor Joseph Ball on occasion of his 
60-th birthday}

\section{Introduction}
\setcounter{equation}{0} 

Let $\cS$ stand for the Schur class of analytic functions mapping the unit
disk ${\mathbb D}$ into $\overline{\mathbb D}$ and let $\cB_\kappa$ be the
set of finite Blaschke products of degree $\kappa$. We denote by 
$\cS_\kappa$ the {\em generalized Schur class} of meromorphic functions of 
the form
\begin{equation}
f(z)=\frac{s(z)}{b(z)},
\label{1.1}
\end{equation}
where $s\in\cS$ and $b\in{\mathcal B}_\kappa$ do not have common zeros (in
particular, $\cS_0=\cS$). Formula \eqref{1.1} is called the Krein-Langer 
representation of a generalized Schur function $f$; the entries $s$ and 
$b$ are determined by $f$ uniquely up to a unimodular constant. 
Via nontangential boundary limits, the $\cS_\kappa$-functions can be identified 
with the functions from the unit ball 
of $L^\infty(\T)$ which admit meromorphic continuation inside the unit disk with 
total pole  multiplicity equal $\kappa$. On the other hand, the $\cS_\kappa$-functions 
$f$ can be characterized as meromorphic functions on $\D$ for which the associated kernel
\begin{equation}
K_f(z,\zeta):=\frac{1-f(z)\overline{f(\zeta)}}{1-z\bar{\zeta}}
\label{1.2}
\end{equation}
has $\kappa$ negative squares on $\rho(f)$, the domain of analyticity of 
$f$: ${\rm sq}_-(K_f)=\kappa$.

\smallskip

The classes ${\cS}_\kappa$ were thoroughly studied in \cite{kl,kl1}, the 
major interpolation results for $\cS_\kappa$-functions can be found in  
\cite{aak,b1,bgr,bh}. The 
objective of this paper is the Nevanlinna-Pick-Carath\'eodory-Fej\'er
interpolation  problem which will be denoted by $\PR$ and
which consists of the following:

\medskip

$\PR$: {\em Given an integer $\kappa\ge 0$,
distinct points $z_1, \ldots, z_k \in
\D$, a tuple ${\bf n}=(n_1,\ldots,n_k)\in{\mathbb N}^k$ and
$|{\bf n}|:=n_1+\ldots+n_k$ complex numbers $f_{i,j}$ $(0\leq j\leq n_i-1; \; 
\; 1\leq i\leq k),$ find all functions $f\in\cS_\kappa$ (if exist) which
are analytic at $z_i$ and satisfy}
\begin{equation}
f^{(j)}(z_i)=j! \, f_{i,j}\qquad (i=1,\ldots,k; \; j=0,\ldots,n_i-1).
\label{1.8}
\end{equation}
Necessary and sufficient conditions for the $\PR$ to have a solution can be given in 
terms of the Pick matrix  of the problem which is determined from interpolation 
data as follows. Let $J_n(z)$ denote the $n\times n$ Jordan block 
with the number $z$ on the main diagonal and let $E_n$ stand for the column vector of
the height $n$ with the first coordinate equals one and other coordinates equal zero:
$$   
J_{n}(z)=\left[\begin{array}{cccc} z&1&&0\\0 &z&\ddots&\\
\vdots&\ddots&\ddots&1\\
 0&\ldots &0&z\end{array}\right],\quad
E_n=\left[\begin{array}{c}1 \\  0 \\ \vdots \\ 0\end{array}\right].
$$
Associated with the tuples ${\bf z}=(z_1,\ldots,z_k)$ and
${\bf n}=(n_1,\ldots,n_k)$ are the matrices
\begin{equation}
T=\left[\begin{array}{ccc}J_{n_1}(z_1) &&
\\ &\ddots & \\ && J_{n_k}(z_k)\end{array}\right],\quad
E=\left[\begin{array}{c} E_{n_1} \\  \vdots \\
E_{n_k}\end{array}\right],
\label{1.5}
\end{equation}
and we arrange the rest of data in the column-vector
\begin{equation}
C=\left[\begin{array}{c} C_1 \\ \vdots \\  C_k\end{array}\right],
\quad\mbox{where}\quad C_i=\left[\begin{array}{c} f_{i,0} \\ \vdots \\
f_{i,n_i-1}\end{array}\right].
\label{1.10}
\end{equation}
Since all the eigenvalues of $T$ fall inside the unit
disk, the Stein equation
\begin{equation}  
P-TPT^*=EE^*-CC^*
\label{1.9}   
\end{equation}
has a unique solution $P$ which is defined via the converging series
\begin{equation}
P=\sum_{j=0}^\infty T^j(EE^*-CC^*)T^{*j}
\label{1.11}  
\end{equation}
and which is called the {\em Pick matrix} of the problem $\PR$. A necessary 
condition for the $\PR$ to have a solution can be obtained as follows.
Given an $f\in\cS_\kappa$, an integer $k\ge 0$ and two
$k$-tuples ${\bf z}=(z_1,\ldots,z_k)\in\rho(f)^k$ and 
${\bf n}=(n_1,\ldots,n_k)\in{\mathbb N}^k$, we define the column-vector
\begin{equation}
M_{\bf n}(f,{\bf z})=\left[\begin{array}{c}
M_{n_1}(f;z_1) \\ \vdots \\ M_{n_k}(f;z_k)\end{array}\right],
\quad\mbox{where}\quad M_{n_i}(f;z_i)=
\left[\begin{array}{c} f(z_i) \\ \frac{f^\prime(z_i)}{1!} \\ \vdots
\\ \frac{f^{(n_i-1)}(z_i)}{(n_i-1)!}\end{array}\right]
\label{2.4a}  
\end{equation}
and the $|{\bf n}|\times|{\bf n}|$ {\em Schwarz-Pick matrix}
\begin{equation}
P_{{\bf n}}(f;{\bf z})=-\frac{1}{4\pi^2}\int_{\Gamma}\int_{\Gamma}
(\xi-T)^{-1}EK_f(\xi,\omega)E^*
(\bar{\omega}-T^*)^{-1} \, d\xi d\bar{\omega},
\label{1.6}
\end{equation}
where $T$ and $E$ are given in \eqref{1.5} and where
$\Gamma\in\D$ is any contour enclosing the points $z_1,\ldots,z_k$
and such that $\overline{{\rm Int} \, \Gamma}\subset \rho(f)$.
Since ${\rm sq_-}(K_f)=\kappa$, the standard approximation arguments show that
the the matrix $P_{{\bf n}}(f;{\bf z})$ has at most $\kappa$ negative eigenvalues:
${\rm sq}_-(P_{{\bf n}}(f;{\bf z}))\le \kappa$. Furthermore, equality 
\begin{equation}
M_{\bf n}(f,{\bf z})=\frac{1}{2\pi i}\int_{\Gamma}(\xi- T)^{-1}Ef(\xi) \,
d\xi
\label{1.7}   
\end{equation}
follows from definitions \eqref{2.4a} and \eqref{1.5} by residue calculus; using this 
equality, one can readily check that the matrix $P_{{\bf n}}(f;{\bf z})$ defined in 
\eqref{1.6} satisfies the Stein equation 
$$
 P_{{\bf n}}(f;{\bf z})-TP_{{\bf n}}(f;{\bf z})T^*=EE^*-M_{\bf 
n}(f;{\bf z})M_{\bf n}(f;{\bf z})^*.
$$
Now we observe that for every solution $f$ of the problem $\PR$, 
the Schwarz-Pick matrix $P_{{\bf n}}(f;{\bf z})$ is equal to $P$,
the Pick matrix of the problem (indeed, if $f$ satisfies interpolation conditions
\eqref{1.8}, it follows from \eqref{1.10} and \eqref{1.7} that $M_{\bf n}(f;{\bf z})=C$; 
thus $P$ and $P_{{\bf n}}(f;{\bf z})$ satisfy the same Stein equation which in turn, has 
a unique solution).  In particular, if $P$ has more than $\kappa$ negative eigenvalues, 
the problem $\PR$ has no solutions. Thus, condition 
$\kappa\geq {\rm sq}_-(P)$ is necessary for the $\PR$ 
to have a solution. On the other hand, if
\begin{equation}
\kappa\ge {\rm sq}_-(P)\quad\mbox{and}\quad\det P\neq 0,
\label{1.12}  
\end{equation}  
then the problem $\PR$ has infinitely many solutions, which
are parametrized by a linear fractional transformation. This is the main result of the 
paper.
\begin{thm}
Let the Pick matrix $P$ of the $\PR$ meet
conditions $(\ref{1.12})$ and let $\Theta=\begin{bmatrix}\Theta_{11} & \Theta_{12}\\
\Theta_{21}& \Theta_{22}\end{bmatrix}$ be the $2\times 2$ rational
matrix-valued function defined by
\begin{equation}
\Theta(z)=I_{2}+(z-1)\left[\begin{array}{c}E^* \\
C^*\end{array}\right](I-zT^*)^{-1}P^{-1}(I-T)^{-1}
\left[\begin{array}{cc}E & -C\end{array}\right].
\label{1.14}
\end{equation}
Then all solutions
$f$ of $\PR$ are parametrized by the 
linear fractional transformation
\begin{equation}
f(z)=\frac{\Theta_{11}(z)S(z)+\Theta_{12}(z)B(z)}
{\Theta_{21}(z)S(z)+\Theta_{22}(z)B(z)},
\label{1.15}
\end{equation}
where the parameters $S\in\cS$ and $B\in{\mathcal B}_{\kappa-{\rm sq}_-(P)}$ 
do not have common zeros and satisfy conditions
\begin{equation}
\Theta_{21}(z_i)S(z_i)+\Theta_{22}(z_i)B(z_i)\neq 0\quad (i=1,\ldots,k).
\label{1.16}  
\end{equation}
Equivalently,
\begin{equation}
f(z)={\bf T}_\Theta[\cE]:=\frac{\Theta_{11}(z)\cE(z)+\Theta_{12}(z)}
{\Theta_{21}(z)\cE(z)+\Theta_{22}(z)},
\label{1.17}  
\end{equation}
where the parameter $\cE\in\cS_{\kappa-{\rm sq}_-(P)}$ satisfies 
\begin{equation}
\Theta_{21}(z_i)\cE(z_i)+\Theta_{22}(z_i)\neq 0
\quad (i=1,\ldots,k)
\label{1.18}
\end{equation}
or has a pole at $z_i$ in case $\Theta_{21}(z_i)\neq\Theta_{22}(z_i)=0$.
\label{T:1.1}
\end{thm}
Equivalence of descriptions \eqref{1.15} and \eqref{1.17} is established 
via the Krein-Langer representation $\cE=\frac{S}{B}$ of the function 
$\cE\in\cS_{\kappa-{\rm sq}_-(P)}$. If $\kappa$ is minimally
possible (i.e., if $\kappa={\rm sq}_-(P)$), then the parameter 
$\cE$ in \eqref{1.17} runs through the Schur class $\cS$; this result can be 
found in \cite{b1,bgr,bh}. A somewhat new point presented here is that 
in case $\kappa>{\rm sq}_-(P)$, some 
solutions of the problem may arise via formula \eqref{1.17} from parameters which 
are not analytic at interpolation nodes. We illustrate 
this possibility by a numerical example.
\begin{ex}
{\rm Let $z_1=0$, $z_2=1/2$, $f_1=1$  
and $f_2=1/2$ so that the interpolation conditions are
\begin{equation}
f(0)=1\quad\mbox{and}\quad f(1/2)=1/2.
\label{6.1}
\end{equation}
The Pick matrix
$P=\begin{bmatrix}0 &&   1/2 \\ 1/2& &  1\end{bmatrix}$ of the problem 
has one negative and one positive eigenvalues and thus, the problem $\PR$
has a solution if and only if $\kappa\ge 1$. Furthermore, substituting
$$
P^{-1}=\left[\begin{array}{rc}-4 & 2\\ 2 & 0\end{array}\right],\quad
T=\begin{bmatrix} 0 & 0 \\ 0 & \frac{1}{2}\end{bmatrix},\quad
E=\begin{bmatrix}1 \\1 \end{bmatrix},\quad C=\begin{bmatrix}1 \\
\frac{1}{2}\end{bmatrix}
$$
into \eqref{1.14} gives
$$
\Theta(z)=\frac{1}{2-z}\begin{bmatrix}3z-2 & 2z(1-z) \\ 2(z-1) & z(3-2z)
\end{bmatrix}.
$$
By Theorem \ref{T:1.1}, all solutions of the problem ${\bf IP}_1$ with interpolation
conditions \eqref{6.1} are parametrized by the formula
\begin{equation}
f(z)=\frac{(3z-2)\cE(z)+2z(1-z)}{2(z-1)\cE(z)+z(3-2z)},
\label{4.1}
\end{equation}  
where $\cE$ belongs to $\cS_0$ and satisfies
\begin{equation}
\cE(0)\neq0\quad\mbox{and}\quad \cE(1/2)\neq 1.
\label{4.2}
\end{equation}
To get all solutions of the problem ${\bf IP}_2$ with interpolation conditions
\eqref{6.1}, we use the formula \eqref{4.1} with parameters $\cE\in\cS_1$ which are
analytic at $z_1=0$ and $z_2=1/2$ and match constraints \eqref{4.2}.
However, since $\Theta_{22}(0)=0$ and $\Theta_{21}(0)=-2\neq 0$, Theorem \ref{T:1.1}
asserts that any function $\cE$ which has a pole at $z_1=0$ and meets the second
constraint in \eqref{4.2}, also leads via \eqref{4.1} to a solution $f$ of the problem.
}
\label{E:1.6}
\end{ex}
Our interest to the ``non-minimal'' problem $\PR$ (where $\kappa>{\rm sq}_-(P)$) 
is motivated by the following reason: if the Pick matrix $P$ of the problem $\PR$
is singular, then the minimally possible $\kappa$ for which the $\PR$ has a solution, 
may be greater than ${\rm sq}_-(P)$.  As we will show in the follow-up paper, 
the description of all solutions for such a degenerate problem can be reduced to a 
family of nondegenerate ``non-minimal'' problems at which point Theorem 
\ref{T:1.1} will come into play. The proof of Theorem \ref{T:1.1} is presented in 
Section 2. In Section 3 we will discuss the divisor-remainder formulation of the 
problem $\PR$ (also considered in \cite{b1,bgr,bh} for $\kappa={\rm sq}_-(P)$).

\section{Proof of Theorem \ref{T:1.1}}
\setcounter{equation}{0}

We first recall some properties of the function $\Theta$ 
defined in \eqref{1.14}. In what follows, $N\{g\}$ stands for the 
total number of zeros of a function $g$ that fall inside $\D$.
\begin{lem}
Let $T$, $E$, $C$ and $P$ be given by \eqref{1.5}, \eqref{1.10} and 
\eqref{1.11}, let $P$ be invertible and let $\Theta(z)$ be defined as in 
\eqref{1.14}. Then
\begin{enumerate}
\item $\Theta(t)$ is $J$-unitary at every point $t\in\T$:
\begin{equation}
\Theta(t)^*J\Theta(t)=\Theta(t)J\Theta(t)^*=J:= \left[\begin{array}{cr}1 & 
0 \\ 0  &-1\end{array}\right].
\label{2.1}   
\end{equation}
\item For $z,\zeta\in\D$,
\begin{equation}
\frac{J-\Theta(z)J\Theta(\zeta)^*}{1-z\bar{\zeta}}=
\left[\begin{array}{c}E^* \\ C^*\end{array}\right]
(I-zT^*)^{-1}P^{-1}(I-\bar{\zeta}T)^{-1}\left[\begin{array}{cc}
E & C\end{array}\right]\label{2.2}
\end{equation}
and
\begin{equation}
\det \, \Theta(z)=\prod_{i=1}^k\left(\frac{(z-z_i)(1-{\bar z}_i)}
{(1-z{\bar z}_i)(1-z_i)}\right)^{n_i}.
\label{2.3}
\end{equation}
\item The function $(zI-T)^{-1}\begin{bmatrix}E & -C\end{bmatrix}\Theta(z)$
is analytic on $\D$.
\item The rational functions $\Theta_{11}$ and $\Theta_{22}$ do not vanish 
on the unit circle and have respectively ${\rm sq}_+(P)$ and ${\rm 
sq}_-(P)$ zeros in $\D$:
\begin{equation}
N\{\Theta_{11}\}={\rm
sq}_+(P)\quad\mbox{and}\quad
N\{\Theta_{22}\}={\rm sq}_-(P).
\label{2.3a}
\end{equation}
\item For every $z\in\C$,
\begin{equation}
|\Theta_{21}(z)|+|\Theta_{22}(z)|>0.
\label{2.3q}
\end{equation}
\item For every $\cE\in\cS_{\widetilde{\kappa}}$, the function 
$f={\bf T}_{\Theta}[\cE]:=\frac{\Theta_{11}\cE+\Theta_{12}}
{\Theta_{21}\cE+\Theta_{22}}$ belongs to $\cS_\kappa$ with 
$\kappa\le \widetilde{\kappa}+{\rm sq}_-(P)$.
\end{enumerate}
\label{L:2.1}  
\end{lem}
{\bf Proof:} Identity \eqref{2.2} follows by a straightforward
calculation (see, e.g., \cite[Section 7.1]{bgr}) based on the identity 
\eqref{1.9}. A similar calculation gives
\begin{eqnarray}
\frac{J-\Theta(\zeta)^*J\Theta(z)}{1-z\bar{\zeta}}&=&
\left[\begin{array}{c}E^* \\ -C^*\end{array}\right](I-T^*)^{-1}P^{-1}
(I-\bar{\zeta}T)^{-1}P(1-zT^*)^{-1}\nonumber\\ 
&&\quad\times P^{-1}(I-T)^{-1}\left[\begin{array}{cc}
E & -C\end{array}\right]\label{2.3c}.
\end{eqnarray}
Identities \eqref{2.1} follow from \eqref{2.2} and \eqref{2.3c}, since 
$\Theta$ is rational and has no poles on $\T$. Equality \eqref{2.3} 
follows from \eqref{1.9} by the standard properties of determinants
(see e.g. \cite[Lemma 2.2]{expar} for the proof).
The third statement of the lemma is yet another consequence 
of identity \eqref{1.9} due to which we have
\begin{eqnarray*}
\begin{bmatrix}E & -C\end{bmatrix}\Theta(z)&=&
\begin{bmatrix}E & -C\end{bmatrix}+(z-1)\begin{bmatrix}E & 
-C\end{bmatrix}\left[\begin{array}{c}E^* \\
C^*\end{array}\right](I-zT^*)^{-1}P^{-1}\\ &&\qquad\qquad\qquad\times(I-T)^{-1}
\left[\begin{array}{cc}E & -C\end{array}\right]\\
&=&  \left(I+(z-1)(P-TPT^*)(I-zT^*)^{-1}P^{-1}(I-T)^{-1}\right)
\left[\begin{array}{cc}E & -C\end{array}\right]\\
&=&(zI-T)P(I-T^*)(I-zT^*)^{-1}P^{-1}(I-T)^{-1}\left[\begin{array}{cc}E
& -C\end{array}\right].
\end{eqnarray*}
Therefore, 
$$
(zI-T)^{-1}\begin{bmatrix}E & -C\end{bmatrix}\Theta(z)=
P(I-T^*)(I-zT^*)^{-1}P^{-1}(I-T)^{-1}\left[\begin{array}{cc}E
& -C\end{array}\right]
$$
and the function on the right hand side is analytic on $\D$.

\smallskip

Equalities \eqref{2.1} imply in particular
\begin{equation}
|\Theta_{11}(t)|^2-|\Theta_{21}(t)|^2=1,\quad
|\Theta_{21}(t)|^2-|\Theta_{22}(t)|^2=-1\quad(t\in\T)
\label{2.3d}
\end{equation}
and thus, $\Theta_{11}$ and $\Theta_{22}$ do not vanish on $\T$. 
For the proof of the second equality in 
\eqref{2.3a}, we refer to  \cite[Theorem 13.2.3]{bgr} or to \cite[Lemma 
4]{Gol}. This equality tells us that if $T$ and $E$ are defined as in 
\eqref{1.5} and $M$ is an arbitrary
vector in $\C^{|{\bf n}|}$ such that the unique solution $R$ (which is 
necessarily Hermitian) of the Stein equation 
\begin{equation}
R-T^*RT=EE^*-MM^*
\label{2.3e}
\end{equation}
 is invertible, then the function 
\begin{equation}
F_M(z)=1-(z-1)M^*(I-zT^*)^{-1}R^{-1}(I-T)^{-1}M
\label{2.3g}
\end{equation}
has ${\rm sq}_-(R)$ zeros inside $\D$:
\begin{equation}
N\{F_M\}={\rm sq}_-(R).
\label{2.3i}
\end{equation}
Let $C$ be the vector associated with 
the problem $\PR$ and decomposed as in \eqref{1.10}. For an 
$\varepsilon>0$, define $C_{\varepsilon}:=C+\varepsilon E$ and the matrix
$P_{\varepsilon}$, a unique solution of the Stein equation 
\begin{equation}
P_{\varepsilon}-TP_{\varepsilon}T^*=EE^*-C_{\varepsilon}C_{\varepsilon}^*.
\label{2.3f}
\end{equation}
Due to the structure \eqref{1.5} of $E$, the above perturbation changes only 
the top entries $f_{i,0}$ in each of the blocks $C_i$ replacing them by 
$f_{i,0}+\varepsilon$. It is clear that there exists $\varepsilon_0$ so 
that for every $\varepsilon\in(0,\varepsilon_0)$,
$$
{\rm sq}_\pm (P_{\varepsilon})={\rm sq}_\pm (P)\quad\mbox{and}\quad 
f_{i,0}+\varepsilon\neq 0\quad\mbox{for}\quad i=1,\ldots,k.
$$
Now we let ${\bf C}_\varepsilon$ to be the block diagonal matrix with
lower triangular Toeplitz diagonal blocks:
$$
{\bf C}_\varepsilon=\left[\begin{array}{ccc} {\bf C}_{\varepsilon,1} & & 
0 \\ & \ddots & \\  0&&  {\bf C}_{\varepsilon,k}
\end{array}\right],\quad
{\bf C}_{\varepsilon,i}=\left[\begin{array}{cccc}f_{i,0}+\varepsilon & 0 & 
\ldots & 0 \\
f_{i,1}& f_{i,0}+\varepsilon & \ddots & \vdots \\ \vdots& \ddots & \ddots 
& 0 \\ f_{i,n_i-1}&  \ldots & f_{i,1} & f_{i,0}+\varepsilon
\end{array}\right].
$$
It is obvious that ${\bf C}_\varepsilon$ is invertible and satisfies 
relations 
$$
{\bf C}_\varepsilon E=C_{\varepsilon}\quad\mbox{and}\quad {\bf 
C}_\varepsilon T=T{\bf C}_\varepsilon.
$$
Multiplying both parts in \eqref{2.3f} by ${\bf C}_\varepsilon^{-1}$ on 
the left and by its adjoint on the right and making use of the two last 
equalities, we get
$$
{\bf C}_\varepsilon^{-1}P_{\varepsilon}{\bf C}_\varepsilon^{-*}
-{\bf 
C}_\varepsilon^{-1}TP_{\varepsilon}{\bf 
C}_\varepsilon^{-*}T^*={\bf C}_\varepsilon^{-1}EE^*{\bf 
C}_\varepsilon^{-*}-EE^*
$$
which can be written in the form \eqref{2.3e} upon setting 
$$
R=-{\bf C}_\varepsilon^{-1}P_{\varepsilon}{\bf C}_\varepsilon^{-*}
\quad\mbox{and}\quad M={\bf C}_\varepsilon^{-1}E.
$$
For this setting, the formula \eqref{2.3g} takes the form
\begin{eqnarray*}
F_\varepsilon(z)&=&1+(z-1)E^*{\bf C}_\varepsilon^{-*}(I-zT^*)^{-1}
({\bf C}_\varepsilon^{-1}P_{\varepsilon}{\bf C}_\varepsilon^{-*})^{-1}
(I-T)^{-1}{\bf C}_\varepsilon^{-1}E\\
&=&1+(z-1)E^*(I-zT^*)^{-1}P_{\varepsilon}^{-1}(I-T)^{-1}E
\end{eqnarray*}
from which it follows that
$$
\lim_{\varepsilon\to 
0}F_\varepsilon(z)=1+(z-1)E^*(I-zT^*)^{-1}P^{-1}(I-T)^{-1}E
=\Theta_{11}(z).
$$
Due to \eqref{2.3i} we have
$$
N\{F_\varepsilon\}={\rm sq}_-(-{\bf 
C}_\varepsilon^{-1}P_{\varepsilon}{\bf C}_\varepsilon^{-*})=
{\rm sq}_+({\bf 
C}_\varepsilon^{-1}P_{\varepsilon}{\bf C}_\varepsilon^{-*})=
{\rm sq}_+(P_{\varepsilon})={\rm sq}_+(P).
$$
Passing to the limit as $\varepsilon\to 0$ implies that $\Theta_{11}$
has ${\rm sq}_+(P)$ zeros in the closed unit disk and since it does not 
have zeros on $\T$, the first equality in \eqref{2.3a} follows.

\smallskip

To prove the fifth statement, note that if 
$\Theta_{21}(z)=\Theta_{22}(z)=0$, then
$\det\Theta(z)=0$ and thus, by formula \eqref{2.3},
inequality \eqref{2.3q} may fail only
at $z\in\{z_1,\ldots,z_d\}$. Let us show that it doesn't. Assuming
that $\Theta_{21}(z_1)=\Theta_{22}(z_1)=0$ we have by \eqref{1.14} and 
\eqref{1.9}
\begin{eqnarray*}  
0&=&\Theta_{21}(z_1)E^*+\Theta_{22}(z_1)C^*\\
&=&C^*+(z_1-1)C^*(I-z_1T^*)^{-1}P^{-1}(I-T)^{-1}
\left(EE^*-CC^*\right)\\
&=&C^*+(z_1-1)C^*(I-z_1T^*)^{-1}P^{-1}(I-T)^{-1}
\left(P-TPT^*\right)\\
&=&C^*(I-z_1T^*)^{-1}P^{-1}(z_1I-T)(I-T)^{-1}P(I-T^*)
\end{eqnarray*}
which is equivalent to
$$
0=C^*(I-z_1T^*)^{-1}P^{-1}(z_1I-T).
$$
Due to the Jordan structure \eqref{2.2} of $T$ and since $z_1\neq z_i$
for $i=2,\ldots,d$, it follows from the last equality that the row-vector
$C^*(I-z_1T^*)^{-1}P^{-1}$ must be of the form
$$
C^*(I-z_1T^*)^{-1}P^{-1}=\begin{bmatrix}\alpha & 0 & \ldots & 0
\end{bmatrix}
$$
and then we have 
\begin{eqnarray*}
\begin{bmatrix}0 & 0 \end{bmatrix}&=&
\begin{bmatrix}\Theta_{21}(z_1) & \Theta_{22}(z_1)\end{bmatrix}\\
&=&\begin{bmatrix}0 & 1 \end{bmatrix}+ (z_1-1)\begin{bmatrix}\alpha & 0 &
\ldots & 0\end{bmatrix}(I-T)^{-1}\left[\begin{array}{cc}
E & -C\end{array}\right]\\
&=& \begin{bmatrix}0 & 1 \end{bmatrix}-\alpha\begin{bmatrix}1 & -f_{1,0}
\end{bmatrix}= \begin{bmatrix}-\alpha & 1+\alpha f_{1,0}\end{bmatrix}
\end{eqnarray*}
which is a contradiction. Thus, inequality \eqref{2.3q} holds for
$z=z_1$ and similarly for $z_2,\ldots,z_d$ which completes the proof of 
the fifth statement.

\smallskip

Finally, let us observe that the function $f={\bf T}_{\Theta}[\cE]$ is 
well defined, i.e., that the function $G_\cE=\Theta_{21}\cE+\Theta_{22}$ 
does not vanish identically for any generalized Schur function $\cE$.
Indeed, if $G_\cE\equiv 0$, then the function 
$-\frac{\Theta_{22}}{\Theta_{21}}=\cE$ belongs to a generalized Schur 
class which is impossible, since due to the second equality in 
\eqref{2.3d}, $\left|\frac{\Theta_{22}(t)}{\Theta_{21}(t)}\right|>1$ at 
every $t\in\T$. Now the last statement in the lemma follows from the 
identity 
$$
G(z)K_f(z,\zeta)G(\zeta)^*=
K_\cE(z,\zeta)+\begin{bmatrix} \cE(\zeta)^* & 1\end{bmatrix}
\frac{J-\Theta(\zeta)^*J\Theta(z)}{1-z\bar{\zeta}}\begin{bmatrix}\cE(z)\\ 1 
\end{bmatrix}
$$
since ${\rm sq}_-\left(\frac{J-\Theta(\zeta)^*J\Theta(z)}{1-z\bar{\zeta}}
\right)\le {\rm sq}_- (P)$, by \eqref{2.3c}. \qed

\medskip

For notational convenience, in what follows
we will often write $f^*$ rather than $\overline{f}$.
\begin{thm}
Let $P$ satisfy conditions \eqref{1.12}, let $\Theta$ be given by \eqref{1.14} and 
let $f$ be a solution of the problem $\PR$.  Then 
\begin{enumerate}
\item The kernel
\begin{equation}
{\bf K}_f(z,\zeta)=\begin{bmatrix} P & 
(I-\bar{\zeta}T)^{-1}\left(E-Cf(\zeta)^*\right)\\
\left(E^*-f(z)C^*\right)(I-zT^*)^{-1}
& K_f(z,\zeta)\end{bmatrix}
\label{2.4}
\end{equation}
has $\kappa$ negative squares on $\rho(f)$.
\item The function $f$ is 
of the form \eqref{1.17} for some $\cE\in\cS_{\kappa-{\rm sq}_- \, (P)}$.
\end{enumerate}
\label{T:2.2}
\end{thm}
{\bf Proof:} Let $\Gamma\in\D$ be any contour enclosing the points $z_1,\ldots,z_k$ 
and such that $\overline{{\rm Int} \, \Gamma}\subset \rho(f)$. Since 
${\rm sq}_-(K_f)=\kappa$, the standard approximation arguments show that 
the kernel 
\begin{equation}
\widetilde{\bf K}_f(z,\zeta)=-\frac{1}{4\pi^2}\int_{\Gamma}\int_{\Gamma}
\begin{bmatrix} (\xi-T)^{-1}E \\ (\xi-z)^{-1}\end{bmatrix} K_f(\xi,\omega)
\begin{bmatrix}E^*(\bar{\omega}-T^*)^{-1} & (\bar{\omega}-\bar{\zeta})^{-1}
\end{bmatrix}d\xi d\bar{\omega}
\label{2.5}
\end{equation}
has at most $\kappa$ negative squares. Since $f$ is a solution of the 
problem $\PR$, we have $P_{{\bf n}}(f;{\bf z})=P$ and $M_{{\bf 
n}}(f;{\bf z})=C$ which together with \eqref{1.6} and \eqref{1.7} lead us to
\begin{equation}
P=-\frac{1}{4\pi^2}\int_{\Gamma}\int_{\Gamma}
(\xi-T)^{-1}EK_f(\xi,\omega)E^*
(\bar{\omega}-T^*)^{-1} \, d\xi d\bar{\omega}
\label{2.6}   
\end{equation}
and 
$$
C=\frac{1}{2\pi i} \int_{\Gamma}(\xi- T)^{-1}Ef(\xi)d\xi.
$$
Using the latter equality along with \eqref{1.2} gives
\begin{eqnarray}
-\frac{1}{4\pi^2}\int_{\Gamma}\int_{\Gamma}(\xi-T)^{-1}EK_f(\xi,\omega)
\frac{d\xi d\bar{\omega}}{\bar{\omega}-\bar{\zeta}}
&=&\frac{1}{2\pi i}\int_{\Gamma}(\xi-T)^{-1}E \, 
\frac{1-f(\xi)f(\zeta)^*}{1-\xi\bar{\zeta}}d\xi\nonumber\\
&=&(I-\bar{\zeta}T)^{-1}E-(I-\bar{\zeta}T)^{-1}Cf(\zeta)^*.\nonumber
\end{eqnarray}
Finally, 
$$
-\frac{1}{4\pi^2}\int_{\Gamma}\int_{\Gamma}
\frac{K_f(\xi,\omega) \, d\xi
d\bar{\omega}}{(\xi-z)(\bar{\omega}-\bar{\zeta})}=K_f(z,\zeta).
$$
Substituting the two last equalities and \eqref{2.6} into \eqref{2.5} and 
comparing the resulting matrix with \eqref{2.4} we conclude that ${\bf 
K}_f(z,\zeta)=\widetilde{\bf K}_f(z,\zeta)$. Therefore, ${\rm sq}_-({\bf K}_f)\le 
\kappa$. 
On the other hand, it follows from \eqref{2.4} that 
${\rm sq}_-({\bf K}_f)\ge {\rm sq}_-(K_f)=\kappa$ which completes the proof of the 
first statement of the theorem. To prove the second, we first note that 
the kernel 
$$
{\bf S}(z,\zeta)=K_f(z,\zeta)-(E^*-f(z)C^*)(I-zT^*)^{-1}
P^{-1}(I-\bar{\zeta}T)^{-1}\left(E-Cf(\zeta)^*\right)
$$
is the Schur complement of the block $P$ in ${\bf K}_f(z,\zeta)$ and 
since ${\rm sq}_-({\bf K}_f)=\kappa$ by the first part, it follows that
\begin{equation}
{\rm sq}_-({\bf S})={\rm sq}_-({\bf K}_f)-{\rm sq}_-(P)=\kappa-{\rm sq}_-(P).
\label{2.8}
\end{equation}
Making use of relation 
$$
K_f(z,\zeta)=\frac{\begin{bmatrix} 1&-f(z)\end{bmatrix}J
\begin{bmatrix}1 \\ -f(\zeta)^*\end{bmatrix}}{1-z\bar{\zeta}}\quad\mbox{where}\quad
J:= \left[\begin{array}{cr}1 &
0 \\ 0  &-1\end{array}\right],
$$
we represent ${\bf S}$ in the form
\begin{eqnarray}
{\bf S}(z,\zeta)&=&\begin{bmatrix}1 & -f(z) 
\end{bmatrix}\left\{\frac{J}{1-z\bar{\zeta}}\right. \label{2.9}\\
&& \quad \left.-\left[\begin{array}{c}E^* \\ C^*\end{array}\right]
(I-zT^*)^{-1}P^{-1}(I-\bar{\zeta}T)^{-1}\left[\begin{array}{cc}
E & C\end{array}\right]\right\}
\begin{bmatrix}1 \\ -f(\zeta)^*\end{bmatrix}\nonumber
\end{eqnarray}
which in turn, can be written as
$$
{\bf S}(z,\zeta)=\begin{bmatrix}1 & -f(z)\end{bmatrix}
\frac{\Theta(z)J\Theta(\zeta)^*}{1-z\bar{\zeta}}\begin{bmatrix}1 \\ 
-f(\zeta)^*\end{bmatrix},
$$
due to \eqref{2.2}. The last equality can be written in terms of the 
pair $\{u,\, v\}$ defined by
\begin{equation}
\begin{bmatrix}v(z) & 
-u(z)\end{bmatrix}=\begin{bmatrix}1 & 
-f(z)\end{bmatrix}\Theta(z)
\label{2.10}
\end{equation}
as 
\begin{equation}
{\bf S}(z,\zeta)=\begin{bmatrix}v(z) & 
-u(z)\end{bmatrix}\frac{J}
{1-z\bar{\zeta}}\begin{bmatrix}v(\zeta)^* \\ -u(\zeta)^*\end{bmatrix}
=\frac{v(z)v(\zeta)^*-u(z)u(\zeta)^*}{1-z\bar{\zeta}}.
\label{2.11}
\end{equation}
Let us show that $v(z)\not\equiv 0$. Indeed, by the first equality in 
\eqref{2.3d}, $\left|\frac{\Theta_{11}(t)}{\Theta_{21}(t)}\right|>1$ for 
every $t\in\T$. Thus, $\frac{\Theta_{11}(t)}{\Theta_{21}(t)}\not\in\cS_\kappa$. 
However, if $v=\Theta_{11}-f\Theta_{21}\equiv 0$ we have
that $f=\frac{\Theta_{11}}{\Theta_{21}}\in\cS_\kappa$
which is a contradiction. Thus, $v\not\equiv 0$ and the function 
$\cE={\displaystyle\frac{u}{v}}$ is meromorphic on $\D$. Equality \eqref{2.11} can be 
written in terms of this function as
$$
{\bf S}(z,\zeta)=v(z)\frac{1-\cE(z)\cE(\zeta)^*}{1-z\bar{\zeta}}v(\zeta)^*=
v(z)K_\cE(z,\zeta)v(\zeta)^*
$$
which together with \eqref{2.8} implies ${\rm sq}_-(K_\cE)={\kappa-{\rm 
sq}_-(P)}$ so that $\cE\in\cS_{\kappa-{\rm sq}_-(P)}$. Finally, it follows from 
\eqref{2.10} that 
$$
f=\frac{{\Theta}_{11}u+{\Theta}_{12}v}{{\Theta}_{21}u+{\Theta}_{22}v}=
\frac{{\Theta}_{11}\cE+{\Theta}_{12}}{{\Theta}_{21}\cE+{\Theta}_{22}}.
$$
Thus, $f$ is of the form \eqref{1.17} with an  $\cE\in\cS_{\kappa-{\rm sq}_-(P)}$ 
which completes the proof of the theorem.\qed

\medskip

Now we will take a closer look at the numerator and the denominator in the 
linear fractional formula \eqref{1.17}. Let 
$\cE$ be a fixed function from $\cS_{\widetilde{\kappa}}$ with the coprime 
factorization
\begin{equation}
\cE(z)=\frac{\num(z)}{\den(z)},\quad \num\in\cS, \; \den\in{\mathcal 
B}_{\widetilde{\kappa}},
\label{2.12}
\end{equation}
and let $\Theta$ be the rational matrix function defined as in 
\eqref{1.14}. Let 
\begin{equation}
U_{S,B}=\Theta_{11}\num+\Theta_{12}\den,\quad
V_{S,B}(z)=\Theta_{21}\num+\Theta_{22}\den,
\label{2.13}
\end{equation}
so that (\ref{1.17}) takes the form 
\begin{equation}
f(z)=\frac{U_{S,B}(z)}{V_{S,B}(z)}.
\label{2.14}
\end{equation}
For the rest of the paper we assume that $V_{S,B}$ has zeros at $z_i$ of 
respective multiplicities $m_i\ge 0$, i.e., that
\begin{equation}
V_{S,B}(z_i)=\ldots=V^{(m_i-1)}_{S,B}(z_i)=0\quad\mbox{and}\quad
V^{(m_i)}_{S,B}(z_i)\neq 0\quad(i=1,\ldots,k).
\label{2.15}
\end{equation}
Since the case where $m_i=0$ is not excluded, the latter assumption is not 
restrictive.
\begin{thm}
Let $P$ be invertible, let $S\in\cS$, $B\in{\mathcal 
B}_{\widetilde{\kappa}}$, let $\Theta$, 
$U_{S,B}$ and $V_{S,B}$ be given as in \eqref{1.14} and \eqref{2.13}.
Then
\begin{enumerate}
\item $N\{V_{S,B}\}={\rm sq}_-(P)+\widetilde{\kappa}$. If in addition,  
$S$ is a finite Blaschke product of degree $m$ (i.e., if $S\in{\mathcal 
B}_m$), then $N\{U_{S,B}\}={\rm sq}_+(P)+m$. 
\item $U_{S,B}$ and $V_{S,B}$ can have a common zero at no point inside
$\D$, but $z_1,\ldots,z_d$.
\item $U_{S,B}$ and $V_{S,B}$ cannot have a common zero at $z_j$ of
multiplicity greater than $n_j$.
\item If $V_{S,B}$ has the zero of multiplicity $m_j>n_j$ at $z_j$,
then $U_{S,B}$ has the zero of multiplicity $n_j$ at $z_j$.
\item If $V_{S,B}$ has the zero of multiplicity $m_j\le n_j$ at $z_j$,
then $U_{S,B}$ has the zero of multiplicity at least $m_j$ at $z_j$.
\end{enumerate}
\label{T:2.3}   
\end{thm}
{\bf Proof:}  By the second equality in \eqref{2.3d},
$|\Theta_{22}(t)|>|\Theta_{21}(t)|$ on $\T$ and since 
$\num\in\cS$ and $\den$ is unimodular on $\T$, it follows that
$$
|\Theta_{22}(t)\den(t)|>|\Theta_{21}(t)\num(t)|
$$
at almost every point $t\in\T$. Then, by Rouch\`{e}'s theorem,
the functions $V_{S,B}=\Theta_{21}\num+\Theta_{22}\den$ and 
$\Theta_{22}\den$ have the same number of zeros in the disk $\{z: \; 
|z|<r\}$ for every $r$ close enough to $1$. Since the rational function 
$\Theta_{22}$ and the finite Blaschke product $\den$ have finitely 
many zeros in $\D$, we let $r\to 1$ to conclude that 
$$
N\{V_{S,B}\}=N\{\Theta_{22}\den\}=N\{\Theta_{22}\}\cdot N\{\den\}=
{\rm sq}_-(P)+\widetilde{\kappa},
$$
where the last equality holds since $N\{\Theta_{22}\}={\rm sq}_-(P)$ (see 
\eqref{2.3a}) and since $N\{\den\}=\widetilde{\kappa}$. Furthermore, 
$|\Theta_{11}(t)|>|\Theta_{12}(t)|$ on $\T$ by the first  equality in 
\eqref{2.3d} and if $S$ is a finite Blaschke product, we have 
$$
|\Theta_{11}(t)\num(t)|>|\Theta_{12}(t)\den(t)|
$$
at every $t\in\T$. Then we use the preceding arguments to conclude
$$
N\{U_{S,B}\}=N\{\Theta_{11}\num\}=N\{\Theta_{11}\}\cdot N\{\num\}= 
{\rm sq}_+(P)+m,
$$
 where the last equality holds since $N\{\Theta_{11}\}={\rm sq}_+(P)$ (see
\eqref{2.3a})  and since $N\{\num\}=m$. This completes the proof of the 
first statement.

\smallskip

To prove the second statement, we write \eqref{2.13} in the matrix form as
\begin{equation}
\left[\begin{array}{c} U_{S,B}(z) \\ V_{S,B}(z)\end{array}\right]=
\Theta(z)\left[\begin{array}{c} \num(z) \\ \den(z)\end{array}\right]
\label{2.19}
\end{equation}   
and assuming that $U_{S,B}(w)=V_{S,B}(w)=0$ at some point $w\in\D$, we get 
$$
\Theta(w)\left[\begin{array}{c} \num(w) \\ \den(w)\end{array}\right]=0
$$
which implies, since $|\num(w)|+|\den(w)|>0$, that
$\det \, \Theta(w)=0$. But by (\ref{2.3}),
$z_1,\ldots,z_k$ are the only zeros of $\det \, \Theta$, which completes 
the proof of the second statement.

\smallskip

Assuming that $U_{S,B}$ and $V_{S,B}$ have the common zero of order
$m_j>n_j$ at $z_j$, we conclude by (\ref{2.19}) that the vector valued
function
$\Theta(z)\left[\begin{array}{c} \num(z) \\ \den(z)\end{array}\right]$
has the zero of multiplicity $m_j>n_j$ at $z_j$. But then, $\det \,
\Theta(z)$ has the zero of multiplicity $m_j>n_j$ at $z_j$, which 
contradicts to  (\ref{2.3}) and completes the proof of the third 
statement.

\smallskip  

By statement (3) in Lemma \ref{L:2.1}, the function   
$$
Q(z):=(zI-T)^{-1}\left[\begin{array}{cc}E & -C\end{array}\right]\Theta(z)
\left[\begin{array}{c} \num(z) \\ \den(z)
\end{array}\right]=(zI-T)^{-1}\left[\begin{array}{cc}E & -C\end{array}\right]
\left[\begin{array}{c} U_{S,B}(z) \\ V_{S,B}(z)\end{array}\right]
$$
(the second equality follows by \eqref{2.19}) is analytic on $\D$ and in 
particular, at $z_1,\ldots,z_k$. 
The block structure \eqref{1.5}, \eqref{1.10} of matrices $T$, $C$ and 
$E$ leads to 
the conformal block structure of $Q$:
\begin{equation}
Q(z)=\left[\begin{array}{c} Q_1(z)\\ \vdots \\  Q_k(z)\end{array}\right],
\quad\mbox{where}\quad
Q_i(z)=(zI-J_{n_i}(z_i))^{-1}\left[E_{n_i}U_{S,B}(z)-C_iV_{S,B}(z)\right]
\label{2.20}
\end{equation}  
and to the conclusion that $Q_i(z)$ is analytic at $z_i$ for
$i=1,\ldots,k$. It is readily seen from the definition of $Q_i$
that the residue of $Q_i$ at $z_i$ is equal to
$$
{\rm Res}_{z=z_i}Q_i(z)=M_{n_i}(U_{S,B};z_i)-
\left[\begin{array}{cccc}f_{i,0} & 0 & \ldots & 0 \\
f_{i,1}& f_{i,0} & \ddots & \vdots \\ \vdots& \ddots & \ddots & 0 \\
f_{i,n_i-1}&  \ldots & f_{i,1} & f_{i,0}
\end{array}\right]M_{n_i}(V_{S,B};z_i)
$$
where $M_{n_i}(U_{S,B};z_i)$ and $M_{n_i}(V_{S,B};z_i)$
are defined in accordance to \eqref{2.4a}.
Since $Q_i$ is analytic at $z_i$ and therefore, ${\rm
Res}_{z=z_i}Q_i(z)=0$, the last displayed equality implies
$$
\left[\begin{array}{cccc}f_{i,0} & 0 & \ldots & 0 \\
f_{i,1}& f_{i,0} & \ddots & \vdots \\ \vdots& \ddots & \ddots & 0 \\
f_{i,n_i-1}&  \ldots & f_{i,1} & f_{i,0}
\end{array}\right]\left[\begin{array}{c}V_{S,B}(z_i) 
\\ \frac{V_{S,B}^{\prime}(z_i)}{1!} \\ \vdots 
\\ \frac{V_{S,B}^{(n_i-1)}(z_i)}{(n_i-1)!}\end{array}\right]=
\left[\begin{array}{c}U_{S,B}(z_i)
\\ \frac{U_{S,B}^{\prime}(z_i)}{1!} \\ \vdots
\\ \frac{U_{S,B}^{(n_i-1)}(z_i)}{(n_i-1)!}\end{array}\right].
$$
Thus, if $m_i\le n_i$, then conditions (\ref{2.15}) force
$$
U_{S,B}^{(j)}(z_i)=0 \quad\mbox{for} \;  j=0,\ldots m_i-1,
$$
which means that $U_{S,B}$ has the zero at $z_i$ of at least the
same multiplicity as $V_{S,B}$ does. If $m_i>n_i$, then the same arguments 
show that
$U_{S,B}$ has zero of multiplicity $\widetilde{m}_i\geq  n_i$ at
$z_i$. If $\widetilde{m}_j> n_i$,
then $z_i$ is a common zero of $U_{S,B}$ and $V_{S,B}$ of multiplicity
greater than $n_j$, which is impossible, by Statement 3. 
Thus, $\widetilde{m}_i= n_i$, which completes the proof of the theorem.\qed

\medskip

{\bf Proof of Theorem \ref{T:1.1}:} By Theorem \ref{T:2.2}, every solution $f$
of the problem $\PR$ is of the form \eqref{1.17}, which is equivalent to representation
\eqref{1.15}. Thus, to prove Theorem \ref{T:1.1}, it suffices to show that a function 
$f$ of the form \eqref{1.15} is a solution of the problem $\PR$ 
problem if and only if the parameters $S\in\cS$ and $B\in{\mathcal B}_\kappa$
meet conditions \eqref{1.16}. To this end, take $f$ in the form 
(\ref{2.14}) and  represent the function $Q_i$ from (\ref{2.20}) as
$$
Q_i(z)=(zI-J_{n_i}(z_i))^{-1}\left[E_{n_i}f(z)-C_i\right]V_{S,B}(z).
$$
If conditions \eqref{1.16} are satisfied, i.e., if $V_{S,B}(z_i)\neq 0$ 
for 
$i=1,\ldots,k$, then $f$ is analytic at $z_1,\ldots,z_d$
and the residue of $Q_i$ at $z_i$ equals
\begin{equation}
0={\rm Res}_{z=z_i}Q_i(z)=
\left[\begin{array}{cccc}r_{i,0} & 0 & \ldots & 0 \\
r_{i,1}& r_{i,0} & \ddots & \vdots \\ \vdots& \ddots & \ddots & 0 \\
r_{i,n_i-1}&  \ldots & r_{i,1} & r_{i,0}
\end{array}\right]\left[\begin{array}{c}V_{S,B}(z_i)
\\ \frac{V_{S,B}^{\prime}(z_i)}{1!} \\ \vdots
\\ \frac{V_{S,B}^{(n_i-1)}(z_i)}{(n_i-1)!}\end{array}\right],
\label{2.21}
\end{equation}
where
\begin{equation}
r_{i,j}=f_{i,j}-\frac{f^{(j)}(z_i)}{j!}\quad (j=0,\ldots,n_i-1).
\label{2.22}
\end{equation}
Note that analyticity of $f$ at $z_i$ is required to establish the second equality 
in \eqref{2.21}; the first holds in any event since $Q_i$ is analytic at 
$z_i$.
Since $V_{S,B}(z_i)\neq 0$, it follows from (\ref{2.21}) that
$$
r_{i,j}=0\quad (j=0,\ldots,n_i-1; \; i=1,\ldots,k),
$$
which is equivalent to (\ref{2.25}), by (\ref{2.22}). Furthermore, 
$V_{S,B}(z)$ has ${\rm sq}_-(P)+(\kappa-{\rm sq}_-(P))=\kappa$ zeros inside $\D$ by 
Theorem \ref{T:2.3} (part (1)) and none of them are canceled by zeros of $U_{S,B}(z)$ 
by statement (2) in the same Theorem \ref{T:2.3} (part (2)). Therefore, $f$ has 
$\kappa$ poles inside $\D$ and it is a generalized Schur function by Lemma 
\ref{L:2.1} (part (4)).  Therefore, $f$ belongs to $\cS_\kappa$ and since it satisfies 
conditions \eqref{2.25}, it solves the problem $\PR$. 

\smallskip 

Let us assume that at least one of the conditions \eqref{1.16} fails, 
i.e., that 
$V_{S,B}(z_i)=0$ for some $i\in\{1,\ldots,k\}$. Then $V_{S,B}(z_i)=0$, by 
statement (5) in Theorem \ref{T:2.3} and after cancellation, it turns out 
that $f$ has
$\kappa^\prime<\kappa$ poles inside $\D$ and therefore, it does not belong to 
$\cS_\kappa$. This is one reason why $f$ is not a solution of the  ${\bf 
P}_\kappa$ 
problem.  Besides, a function $f\in\cS_{\kappa^\prime}$ cannot satisfy all 
the interpolation  conditions in \eqref{2.25}. If it had, then by virtue 
of Theorem \ref{T:2.2} it 
would be of the form $f={\bf T}_\Theta[\cE^\prime]$ for some 
$\cE^\prime\in\cS_{\kappa^\prime-{\rm sq}_-(P)}$ and the same coefficient matrix 
$\Theta$. Since the map $\cE\to {\bf T}_\Theta[\cE]$ is invertible, we would have 
$\cE\equiv\cE^\prime$ which is impossible since the latter functions have different 
numbers of poles in $\D$. This completes the proof of Theorem 
\ref{T:1.1}.\qed

\medskip

In conclusion we consider the functions $f$ obtained via the formula 
\eqref{1.15} from the parameters $\{S, \, B\}$ which fail to satisfy all the 
conditions \eqref{1.16}. We will be interested in two questions: 
how many negative squares $f$ may lose and which 
interpolation  conditions it still satisfies. In addition to the 
tuple ${\bf n}=(n_1,\ldots,n_k)\in{\mathbb Z}_+^k$ associated with the 
problem $\PR$ we consider another tuple ${\bf 
m}=(m_1,\ldots,m_k)\in{\mathbb Z}_+^k$ and introduce
\begin{equation}
{\mathcal I}_+:=\{i\in\{1,\ldots,k\}: \; \; n_i> m_i\},\quad
{\mathcal I}_-:=\{i\in\{1,\ldots,k\}: \; \; n_i< m_i\},
\label{3.6}
\end{equation}
\begin{equation}
{\mathcal I}_0=\{i\in\{1,\ldots,k\}: \; \; n_i= m_i\},\quad
\gamma_{\bf m}:={\displaystyle\sum_{i=1}^k\min\{m_i, \; n_i\}}.
\label{3.7}
\end{equation}
\begin{thm}
Let $f$ be of the form \eqref{1.15} with $S\in\cS$ and $B\in{\mathcal 
B}_{\widetilde{\kappa}}$ having no common zeros on $\D$ and such that 
\begin{equation}
(\Theta_{21}S+\Theta_{22}B)^{(j)}(z_i)=0 \quad (i=1,\ldots,k; \;
j=0,\ldots,m_i-1)
\label{2.23}
\end{equation}
and 
\begin{equation}
(\Theta_{21}S+\Theta_{22}B)^{(m_i)}(z_i)\neq 0 \quad (i=1,\ldots,k).
\label{2.24}
\end{equation}
Then $f$ belongs to the class $\cS_{\widetilde{\kappa}+{\rm sq}_-(P)-\gamma_{\bf 
m}}$, where $\gamma_{\bf m}$ is given in 
\eqref{3.7}. Furthermore, $f$ has a pole  of multiplicity $m_i-n_i$
at $z_i$ if $i\in{\mathcal I}_-$, and satisfies interpolation conditions
\begin{equation}
f^{(j)}(z_i)=j! \, f_{i,j} \quad (i\in{\mathcal I}_+; \; j=0,\ldots,n_i-m_i-1).
\label{2.25}  
\end{equation}
\label{T:2.4}
\end{thm}
{\bf Proof:} We take $f$ in the form \eqref{2.14} with $U_{S,B}$ and 
$V_{S,B}$ defined as 
in \eqref{2.13}. Conditions \eqref{2.23}, \eqref{2.24} say that $V_{S,B}$ 
has zeros of
order $m_i$ at $z_i$ for $j=1,\ldots,k$. If $m_i\le n_i$, then 
$U_{S,B}$ has zero of 
order at least $m_i$ at $z_i$ (statement (5) in Theorem \ref{T:2.3}) and 
therefore 
$f$  admits an analytic continuation to $z_i$. If $m_i>n_i$, then the 
same arguments show that $U_{S,B}$ has zero of multiplicity $n_i$ at $z_i$ 
(statement (4) in Theorem  \ref{T:2.3}) and after cancellation $f$ will 
have a pole of multiplicity $m_i-n_i$ at $z_i$. By statement (1) in 
Theorem \ref{T:2.3}, the total number of  zeros of $V_{S,B}$ inside $\D$ 
is  equal to $\widetilde{\kappa}+{\rm sq}_-(P)$. Therefore  
$\widetilde{\kappa}+{\rm sq}_-(P)-|{\bf m}|$ zeros fall into   
$\D\setminus\{z_1,\ldots,z_k\}$ and cannot be canceled by zeros of 
$U_{S,B}$ by  statement (2) in Theorem \ref{T:2.3}. After all zero cancellations, 
the function $V_{S,B}$ will have $m_i-n_i$ zeros every $z_i$ for 
$i\in{\mathcal I}_-$ and still $\widetilde{\kappa}+{\rm sq}_-(P)-|{\bf m}|$ 
zeros in $\D\setminus\{z_1,\ldots,z_k\}$. Thus, the function  
$f=\frac{U_{S,B}}{V_{S,B}}$ will have
$$
\widetilde{\kappa}+{\rm sq}_-(P)-|{\bf m}|+\sum_{i\in{\mathcal I}_-}(m_i-n_i)
=\widetilde{\kappa}+{\rm sq}_-(P)-\gamma_{\bf m}
$$
poles inside $\D$. By statement (4) in Lemma \ref{L:2.1}, $f$ is a 
generalized Schur function and therefore, it belongs to $\cS_{\widetilde{\kappa}+{\rm 
sq}_-(P)-\gamma_{\bf m}}$. Furthermore, if $i\in{\mathcal I}_+$, then $f$ 
is analytic at $z_i$ and therefore equality \eqref{2.21} holds. 
Since $V_{S,B}^{(m_i)}(z_i)\neq 0$ by \eqref{2.24}, it follows from 
(\ref{2.21}) that
$$
r_{i,j}=0\quad (i\in{\mathcal I}_+; \; j=0,\ldots,n_i-m_i-1),
$$
which is equivalent to (\ref{2.25}), by (\ref{2.22}).

\section{The divisor-remainder version}
\setcounter{equation}{0}

The problem $\PR$ can be formulated in the divisor-remainder form \eqref{3.4}
as follows. Let $H^\infty_\kappa$ be the set of all functions $f$ of the form  
\eqref{1.1} where $s\in H^\infty$ and $b\in{\mathcal B}_\kappa$ may have 
common zeros. From this definition it follows that 
$\cS_\kappa=(H^\infty_\kappa\backslash 
H^\infty_{\kappa-1})\cap {\mathcal B}L^\infty$ where ${\mathcal B}L^\infty$
denotes the unit ball of $L^\infty(\T)$. Let 
$\varphi\in H^\infty$ be any function satisfying interpolation 
conditions \eqref{1.8}:
\begin{equation}
\varphi^{(j)}(z_i)=j! \, f_{i,j}\qquad (i=1,\ldots,k; \;
j=0,\ldots,n_i-1),
\label{3.2}
\end{equation}
and let $\theta$ be a finite Blaschke product defined by
\begin{equation}
\theta(z)=\prod_{i=1}^{k}\left(\frac{z-z_i}{1-z\bar{z}_i}
\right)^{n_i}.
\label{3.3}     
\end{equation}
\begin{prop}
A function $f$ is a solution of the problem $\PR$ if and only if it belongs to 
$\cS_\kappa$ and admits a representation
\begin{equation}
f(z)=\varphi(z)+\theta(z)h(z)\quad \mbox{for some}\quad h\in
H_\kappa^\infty.
\label{3.4}  
\end{equation}
\label{P:3.1}
\end{prop}
{\bf Proof:} If $f=s/b$ (where $s\in\cS$ and $b\in{\mathcal B}_\kappa$) belongs 
to $\cS_\kappa$ and satisfies conditions \eqref{1.8}, then by \eqref{3.2}, the 
function $s-\varphi b$ belongs to $H^\infty$ and satisfies the
homogeneous conditions
\begin{equation}
(s-\varphi b)^{(j)}(z_i)=0\qquad (i=1,\ldots,k; \; j=0,\ldots,n_i-1)
\label{3.5}
\end{equation}
By the maximum modulus principle, the function
$(s-\varphi b)/\theta$ belongs to $H^\infty$ where
$\theta$ is
defined in \eqref{3.3}. Since $b\in{\mathcal B}_\kappa$, the function
$h:={\displaystyle\frac{s-\varphi b}{\theta b}}$ belongs to
$H^\infty_\kappa$ and therefore $f$ can be represented as in \eqref{3.4}, 
since   
$$
f=\frac{s}{b}=\varphi+\theta\cdot \frac{s-\varphi b}{\theta
b}=\varphi+\theta h.
$$
Conversely, let $f\in\cS_\kappa$ be of the form \eqref{3.4}. Since $f$
has $\kappa$ poles and $\varphi\in
H^\infty$, it follows that $h$ has $\kappa$ poles in $\D$ none of which
are in $\{z_1,\ldots,z_k\}$, the zero set of $\theta$. Therefore, $h$ is
analytic at $z_1,\ldots,z_k$. Therefore, the function $f-\varphi=\theta h$
satisfies the homogeneous conditions \eqref{3.5}, so that $f$ satisfies
\eqref{1.8} due to \eqref{3.5}.\qed

\medskip

By Proposition \ref{P:3.1}, the solution set for the problem $\PR$ is 
equal to $(\varphi+\theta H^\infty_\kappa)\cap \cS_\kappa$. Thus if the Pick 
matrix $P$ of the problem meets conditions \eqref{1.12}, then the set 
$(\varphi+\theta H^\infty_\kappa)\cap \cS_\kappa$ is not empty (and is 
parametrized as in Theorem \ref{T:1.1}) and therefore, a larger set 
\begin{equation}
\Omega_\kappa(\varphi,\theta):=(\varphi+\theta H^\infty_\kappa)\cap 
{\mathcal B}L^\infty=(\varphi+\theta H^\infty_\kappa)\cap 
\left(\bigcup_{\alpha=0}^\kappa 
\cS_\alpha\right)
\label{3.6a}
\end{equation}
is not empty. The second equality in \eqref{3.6a} is easily verified: since 
$\varphi,\theta\in 
H^\infty$, it follows that $\varphi+\theta H^\infty_\kappa\subset H^\infty_\kappa$ and 
on the 
other hand, $H^\infty_\kappa\cap {\mathcal B}L^\infty=\cup_{\alpha=0}^\kappa \cS_\alpha$.
Clearly, the elements of $\Omega_\kappa(\varphi,\theta)$ are 
solutions of certain $L^\infty$-norm constraint interpolation problem; the 
next theorem characterizes $\Omega_\kappa(\varphi,\theta)$ in terms of 
the kernel ${\bf K}_f$ defined in \eqref{2.4} as well as in terms of 
the linear fractional transformation ${\bf T}_\Theta$ defined in  \eqref{1.17}.
\begin{thm}
Let the Pick matrix $P$ defined in \eqref{1.11} be invertible and let 
$\kappa\ge {\rm sq}_-(P)$. Let $\Theta$ be given by \eqref{1.14}, let 
$\varphi$ be an $H^\infty$-function satisfying conditions \eqref{3.2}, let 
$\theta$ be given by \eqref{3.3} and let $f$ be a function meromorphic on $\D$.
The following are equivalent:
\begin{enumerate}
\item $f$ belongs $(\varphi+\theta H^\infty_\kappa)\cap
{\mathcal B}L^\infty$.
\item ${\rm sq}_-({\bf K}_f(z,\zeta))\le \kappa$ where the kernel
${\bf K}_f$ is defined in \eqref{2.4}.
\item $f={\bf T}_\Theta[\cE]$ for some $\cE\in H^\infty_{\kappa-{\rm
sq}_-(P)}\cap {\mathcal B}L^\infty$.
\end{enumerate}
\label{T:3.1}
\end{thm}
{\bf Proof:} ${\bf (1)\Rightarrow(2)}$. Let us assume that $f$ is of the form 
\eqref{3.4} and belongs to $\cS_{\widetilde{\kappa}}$ 
for some $\widetilde{\kappa}\le \kappa$. Let us assume that the function 
$h$ in representation \eqref{3.4} has poles of multiplicities $m_i$ at
$z_i$ for $i=1,\ldots,k$ (the case where $m_i=0$ is not excluded). Let 
${\mathcal I}_\pm$, ${\mathcal I}_0$ and $\gamma_{\bf m}$ be defined as in 
\eqref{3.6}, \eqref{3.7}. Since $h\in H^\infty_\kappa$, it may have at most 
$\kappa-|{\bf m}|$ 
poles outside the set $\{z_1,\ldots,z_k\}$. After cancellation of the 
poles of $h$ with the zeros of $\theta$, we obtain the following 
representation for $f$:
\begin{equation}
f(z)=\varphi(z)+\widetilde{\theta}(z)\widetilde{h}(z),
\label{3.8}
\end{equation}
where 
\begin{equation}
\widetilde{\theta}(z)=\prod_{i\in{\mathcal I}^+}\left(\frac{z-z_i}
{1-z\bar{z}_i}\right)^{n_i-m_i}\quad\mbox{and}\quad
\widetilde{h}(z)=h(z)\cdot \prod_{i\in{\mathcal I}^-}\left(\frac{z-z_i}
{1-z\bar{z}_i}\right)^{m_i-n_i}.  
\label{3.8b}
\end{equation}
It is clear that $\widetilde{h}$ has poles of multiplicities $m_i-n_i$ at
$z_i$ for every $i\in {\mathcal I}_-$ and at most $\kappa-|{\bf 
m}|$ poles in $\D\setminus\{z_1,\ldots,z_k\}$. It follows from \eqref{3.7}
that $f$ has the same poles and since $f$ belongs to 
$\cS_{\widetilde{\kappa}}$, we get  
\begin{equation}
\widetilde{\kappa}\le \sum_{i\in{\mathcal I}_-}(m_i-n_i)+\kappa-|{\bf
m}|=\kappa-\gamma_{\bf m}.
\label{3.8a}
\end{equation}
On the other hand, $\widetilde{h}$ is analytic at every $z_i$ for 
$i\in{\mathcal I}_+$ and therefore, $f$ of the form \eqref{3.8} 
satisfies interpolation conditions
\begin{equation}
f^{(j)}(z_i)=\varphi^{(j)}(z_i)=j! \, f_{i,j}\qquad (i\in{\mathcal I}_+; \; 
j=0,\ldots,n_i-m_i-1).
\label{3.9}
\end{equation}
Therefore, $f$ is a solution of the problem ${\bf 
IP}_{\widetilde{\kappa}}$ with interpolation conditions \eqref{3.9}.
The Pick matrix $\widetilde{P}$ of this interpolation problem is a 
principal submatrix of the Pick matrix $P$ of the original $\PR$, and for 
a suitable permutation matrix $U$, we have
\begin{equation}
UPU=\begin{bmatrix}P_1 & P_2^* \\ P_2 & 
\widetilde{P}\end{bmatrix}.
\label{3.10}   
\end{equation}
Furthermore, associating the matrices $\widetilde{T}$, $\widetilde{E}$
and  $\widetilde{C}$ with the problem ${\bf IP}_{\widetilde{\kappa}}$
via formulas \eqref{1.5} and \eqref{1.11}, it is easy to check the 
block decompositions
\begin{equation}
UTU=\begin{bmatrix}T_1& 0 \\ T_2 &\widetilde{T}\end{bmatrix},\quad
UE=\begin{bmatrix}E_1 \\ \widetilde{E}\end{bmatrix},\quad
UE=\begin{bmatrix}C_1 \\ \widetilde{C}\end{bmatrix}
\label{3.11}
\end{equation}
conformal with \eqref{3.10}. By Theorem \ref{T:2.2} applied to the problem 
${\bf IP}_{\widetilde{\kappa}}$, the kernel 
\begin{equation}
\widetilde{\bf K}_f(z,\zeta)=\begin{bmatrix} \widetilde{P} &
(I-\bar{\zeta}\widetilde{T})^{-1}\left(\widetilde{E}-\widetilde{C}
f(\zeta)^*\right)\\
\left(\widetilde{E}^*-f(z)\widetilde{C}^*\right)(I-z\widetilde{T}^*)^{-1}
& K_f(z,\zeta)\end{bmatrix}
\label{3.12}   
\end{equation}
has $\widetilde{\kappa}$ negative squares on $\rho(f)$. Due to 
\eqref{3.10}--\eqref{3.12}, the kernel ${\bf K}_f$ defined in \eqref{2.4}
can be represented as 
\begin{equation}
\begin{bmatrix}U & 0 \\ 0 & 1\end{bmatrix}
{\bf K}_f(z,\zeta)\begin{bmatrix}U^* & 0 \\ 0 & 1\end{bmatrix}
=\begin{bmatrix}P_1 & B(\zeta)^* \\ B(z) & \widetilde{\bf K}_f(z,\zeta)
\end{bmatrix},
\label{3.13}  
\end{equation}
where $B(z)$ is analytic on $\rho(f)$ (the explicit formula for $B$ is not 
that important for now). Then the number of negative squares of the 
kernel on the right hand side of \eqref{3.13} can be estimated as follows
\begin{equation}
{\rm sq}_-\left(\begin{bmatrix}P_1 & B(\zeta)^* \\ B(z) & 
\widetilde{\bf K}_f(z,\zeta)
\end{bmatrix}\right) \le {\rm sq}_-(\widetilde{\bf 
K}_f)+d=\widetilde{\kappa}+d,
\label{3.14}
\end{equation}
where $d$ is the size of the square matrix $P_1$ (see \cite[Proposition 
4.1]{bkh1} for the proof). The number $d$ is equal to the difference 
between the sizes of the Pick matrices $P$ and $\widetilde{P}$ or, which 
is the same, to the difference between the numbers of interpolation 
conditions in \eqref{1.8} and \eqref{3.9}. Thus,
\begin{equation}
d=|{\bf n}|-\sum_{i\in{\mathcal I}_+}(n_i-m_i)=\sum_{i\in{\mathcal 
I}_-\cup{\mathcal I}_0}n_i+\sum_{i\in{\mathcal I}_+}m_i=\gamma_{\bf m}
\label{3.15}  
\end{equation}
Now we combine \eqref{3.8a}, \eqref{3.14} and \eqref{3.15} to conclude 
from \eqref{3.13} that
$$
{\rm sq}_-({\bf K}_f)={\rm sq}_-\left(\begin{bmatrix}P_1 & B(\zeta)^* \\ 
B(z) & \widetilde{\bf K}_f(z,\zeta)
\end{bmatrix}\right) \le 
\widetilde{\kappa}+d=\widetilde{\kappa}+\gamma_{\bf m}\le 
\kappa-\gamma_{\bf m}+\gamma_{\bf m}=\kappa,
$$
which completes the proof of the implication ${\bf (1)\Rightarrow(2)}$.

\smallskip

${\bf (2)\Rightarrow(3)}$. Let us assume that ${\rm sq}_-({\bf 
K}_f)=\widetilde{\kappa}\le \kappa$. Then by the arguments used in the 
proof of Theorem \ref{T:2.2} to derive statement (2) from statement (1)
we conclude that $f$ is of the form $f={\bf T}_\Theta[\cE]$ for some 
$\cE\in\cS_{\widetilde{\kappa}-{\rm sq}_-(P)}$. Since 
$\cS_{\widetilde{\kappa}-{\rm sq}_-(P)}\subset H^\infty_{\widetilde{\kappa}-{\rm
sq}_-(P)}\cap {\mathcal B}L^\infty\subset H^\infty_{\kappa-{\rm
sq}_-(P)}\cap {\mathcal B}L^\infty$, the proof is completed.

\smallskip

${\bf (3)\Rightarrow(1)}$. Let $f$ be of the form $f={\bf T}_\Theta[\cE]$ for some 
$\cE\in\cS_{\widetilde{\kappa}}$ where $\widetilde{\kappa}\le \kappa-{\rm sq}_-(P)$.
Then we equivalently can take $f$ in the form \eqref{1.15} with $S\in\cS$ and 
$B\in{\mathcal B}_{\widetilde{\kappa}}$ having no common zeros on $\D$. Let 
$m_1,\ldots,m_k$ be the integers uniquely determined from conditions 
\eqref{2.23} and \eqref{2.24}. Then we conclude by Theorem \ref{T:2.4} that 
$f$ satisfies interpolation conditions \eqref{2.25} (or \eqref{3.9} which is the same)
and belongs to the class $\cS_{\kappa_1}$ where 
\begin{equation}
\kappa_1={\widetilde{\kappa}+{\rm sq}_-(P)-\gamma_{\bf m}}.
\label{3.16}
\end{equation}
Thus, $f$ solves the problem ${\bf IP}_{\kappa_1}$ with interpolation conditions \eqref{3.9} and 
therefore by virtue of Proposition \ref{P:3.1}, it admits a representation 
\eqref{3.8} with 
$\widetilde{\theta}(z)$ defined as in \eqref{3.8b} and some $\widetilde{h}\in 
H^\infty_{\kappa_1}$. From \eqref{3.3} and \eqref{3.8b} we observe that 
the ratio $\theta_1:=\theta/\widetilde{\theta}$ is a finite Blaschke product of degree
$$
\deg \, \theta_1=\sum_{i\in{\mathcal I}_+}m_i+\sum_{i\in{\mathcal I}_-\cup {\mathcal 
I}_0}n_i=\gamma_{\bf m}
$$
and therefore the function $h:=\widetilde{h}/\theta_1$ belongs to $H^\infty_{\kappa_2}$ 
where $\kappa_2=\kappa_1+\gamma_{\bf m}$ and due to \eqref{3.16}, we have
$\kappa_2=\widetilde{\kappa}+{\rm sq}_-(P)\le \kappa$,
so that $h\in H^\infty_{\kappa}$. Now we get from \eqref{3.8}
$$
f=\varphi+\widetilde{\theta}\cdot\widetilde{h}=\varphi+\theta\cdot 
\frac{\widetilde{h}}{\theta_1}=\varphi+\theta h.
$$
Since $f\in\cS_{\kappa^\prime}$ and $h\in H^\infty_{\kappa}$, it follows that 
$$
f\in (\varphi+\theta H^\infty_\kappa)\cap \cS_{\kappa^\prime}\subset (\varphi+\theta 
H^\infty_\kappa)\cap{\mathcal B}L^\infty
$$
which completes the proof of the theorem.\qed 

\medskip

As corollary, we obtain the following ``if and only if'' version of Theorem 
\ref{T:2.2}.
\begin{thm}
Let $P$ satisfy conditions \eqref{1.12}, let $\Theta$ be given by \eqref{1.14} and
let $f\in\cS_\kappa$. The following are equivalent:
\begin{enumerate}
\item $f$ is a solution of the problem $\PR$.
\item ${\rm sq}_-({\bf K}_f(z,\zeta))=\kappa$ where the kernel
${\bf K}_f$ is defined in \eqref{2.4}.
\item $f={\bf T}_\Theta[\cE]$ for some $\cE\in \cS_{\kappa-{\rm sq}_- \, (P)}$.
\end{enumerate}
\label{T:3.5}  
\end{thm}
{\bf Proof:} Implications  ${\bf (1)\Rightarrow(2)\Rightarrow(3)}$ are proved in 
Theorem \ref{T:2.2}. If $f={\bf T}_\Theta[\cE]$ for some $\cE\in \cS_{\kappa-{\rm 
sq}_- \, (P)}$, then $f\in (\varphi+\theta H^\infty_\kappa)\cap
{\mathcal B}L^\infty$, by Theorem \ref{T:3.1}, where $\phi$ and $\theta$ are the 
functions associated with the problem $\PR$. By the assumption of the theorem,
$f\in\cS_\kappa$ and therefore 
$$
f\in (\varphi+\theta H^\infty_\kappa)\cap
{\mathcal B}L^\infty\cap \cS_\kappa=(\varphi+\theta H^\infty_\kappa)\cap  
\cS_\kappa
$$
and the latter set coincides with the solution set for the problem $\PR$, by 
Proposition \ref{P:3.1}. This completes the proof of the implication ${\bf 
(3)\Rightarrow(1)}$ and therefore, of the theorem.\qed

\bibliographystyle{amsplain}
\providecommand{\bysame}{\leavevmode\hbox to3em{\hrulefill}\thinspace}

\end{document}